\documentclass[11pt,a4paper]{article}

\usepackage[T1]{fontenc}
\usepackage[utf8]{inputenc}
\usepackage{lmodern}
\usepackage[a4paper,margin=28mm]{geometry}
\usepackage{amsmath,amssymb,amsthm,mathtools}
\usepackage{microtype}
\usepackage{graphicx,tikz}
\usepackage{needspace,placeins}
\usepackage[hidelinks]{hyperref}
\newtheorem{theorem}{Theorem}[section]
\newtheorem{lemma}[theorem]{Lemma}
\newtheorem{proposition}[theorem]{Proposition}
\newtheorem{corollary}[theorem]{Corollary}
\theoremstyle{definition}
\newtheorem{example}[theorem]{Example}

\numberwithin{equation}{section}
\hypersetup{
 pdftitle={Compatible zero-sum partitions and applications to group distance magic labeling},
 pdfauthor={Ahmet Batal and Sylwia Cichacz},
 pdfsubject={Compatible zero-sum partitions of finite Abelian groups and their applications}
}

\makeatletter
\renewcommand{\@maketitle}{%
  \begin{center}%
    {\LARGE\@title\par}%
    \ifx\@author\@empty\else
      \vskip 1em {\normalsize\@author\par}%
    \fi
    \ifx\@date\@empty\else
      \vskip .5em {\normalsize\@date\par}%
    \fi
  \end{center}%
  \vskip .7em
}
\makeatother

\title{Compatible zero-sum partitions of finite Abelian groups\\and their applications}
\author{{ Ahmet Batal$^{1,2}$\thanks{Ahmet Batal's research was partially supported by program ``Excellence initiative -- research university'' (IDUB) for the AGH University.}, {Sylwia Cichacz$^{2}$\thanks{Sylwia Cichacz's research was supported by  the AGH University of Krakow under grant no. 16.16.420.054, funded by the Polish Ministry of Science and Higher Education.}}}
	\footnote{e-mails: cichacz@agh.edu.pl, ahmetbatal@iyte.edu.tr}\\
\normalsize $^1$Department of Mathematics, Izmir Institute of Technology, 35430, Urla, Izmir, Turkey\\
\normalsize $^2$AGH University of Krakow, Poland \\
}

\begin{document}
\maketitle

\begin{abstract}
Let $(\Gamma,+)$ be a finite Abelian group and let
$g\in\Gamma\setminus\{0\}$. A partition of $\Gamma$ is called
$g$-compatible if no part contains both $x$ and $x+g$.
For a prescribed involution $h$ and odd $r\ge3$, we prove that
$\Gamma$ admits an $h$-compatible partition into zero-sum
$r$-element sets if and only if $r$ divides $|\Gamma|$ and
$\Gamma$ has at least three involutions. Moreover, the partition
can be chosen so that the multigraph whose vertices are its
parts and whose edges are the pairs $\{x,x+h\}$ is Hamiltonian
and has a proper edge coloring with $r$ colors.
We also give examples, for groups of odd order and for even
block sizes, in which zero-sum partitions exist but compatible
ones do not.

We give two applications of these partitions. For block size
three, we obtain complete mappings whose cycles are precisely
the parts of the partition. We also construct Hamiltonian
$r$-regular $\Gamma$-distance-magic graphs admitting $0$ and
every involution of $\Gamma$ as magic constants.
\end{abstract}



\section{Introduction}

Let $(\Gamma,+)$ be a finite Abelian group of order $n$, with
identity element $0$. A nonzero element $h\in\Gamma$ satisfying
$2h=0$ is called an \emph{involution}. We denote the set of
involutions of $\Gamma$ by $I(\Gamma)$. The sum of all elements
of $\Gamma$ is determined by its involutions. We recall the
following identity.

\begin{lemma}[{\cite[Lemma~8]{ref_ComNelPal}}]\label{CN}
Let $\Gamma$ be a finite Abelian group. Then
\begin{equation}\label{eq:groupsum}
 \sum_{x\in\Gamma}x=
 \begin{cases}
 \iota,&I(\Gamma)=\{\iota\},\\
 0,&|I(\Gamma)|\ne1.
 \end{cases}
\end{equation}
\end{lemma}

A subset $S\subseteq\Gamma$ is called \emph{zero-sum} if
$\sum_{x\in S}x=0$. A \emph{zero-sum partition of block size $r$}
is a partition $\mathcal P$ of $\Gamma$ into $r$-element sets such
that
\[
 \sum_{x\in T}x=0
\]
for every $T\in\mathcal P$.
For $r\ge3$, Cichacz~\cite[Theorem~12]{Cichacz2018} proved that
such a partition exists if and only if $r$ divides $|\Gamma|$
and $|I(\Gamma)|\ne1$.

For a nonzero element $g\in\Gamma$, we call a partition
$\mathcal P$ of $\Gamma$ \emph{$g$-compatible} if
\begin{equation}\label{eq:h-compatible}
 T\cap(T+g)=\varnothing
\end{equation}
for every $T\in\mathcal P$, where $T+g=\{x+g:x\in T\}$.
Thus $x$ and $x+g$ must belong to different parts.

The conditions for an ordinary zero-sum partition need not be
sufficient for a $g$-compatible one, even when $\Gamma$ has odd
order. Indeed, let
\[
 \Gamma=\mathbb Z_3\oplus\mathbb Z_5,\qquad
 g=(1,0),\qquad r=3.
\]
Any two distinct elements with the same second coordinate differ
by $g$ or $-g$. Hence the second coordinates of a $g$-compatible
triple are distinct. If the triple is zero-sum, these coordinates
must form one of the sets
\[
 \{0,1,4\}\qquad\text{or}\qquad\{0,2,3\},
\]
which are the only three-element zero-sum subsets of
$\mathbb Z_5$. Thus every part of a $g$-compatible zero-sum
partition into triples would contain an element of
$\mathbb Z_3\oplus\{0\}$. Such a partition would have five
parts, whereas this subgroup has only three elements.
Consequently, no such partition exists.

For odd block sizes, we prove that the conditions for an ordinary
zero-sum partition are sufficient for $h$-compatibility with every
prescribed involution $h$. The existence of $h$ restricts us to
groups of even order. Moreover,
Lemma~\ref{CN} excludes groups with exactly one involution,
since a zero-sum partition requires
$\sum_{x\in\Gamma}x=0$. Therefore $|I(\Gamma)|>1$, or
equivalently $|I(\Gamma)|\ge3$, is a necessary condition for
the partitions considered here.

Even under this condition, even block sizes give further
obstructions. For any even $r\ge4$, take
\[
 \Gamma=\mathbb Z_2\oplus\mathbb Z_r,\qquad h=(1,0).
\]
Then $r\mid|\Gamma|$ and $|I(\Gamma)|=3$.
Every $h$-compatible $r$-element subset contains exactly one
element with each second coordinate $0,1,\ldots,r-1$. The sum
of these coordinates satisfies
\[
 \sum_{j=0}^{r-1}j=\frac{r(r-1)}2
 \equiv\frac r2\not\equiv0\pmod r.
\]
Thus no such subset is zero-sum. This shows that the restriction
to odd block sizes is necessary for the characterization below.
Throughout the remainder of the paper, $r\ge3$ is odd, and
$h\in I(\Gamma)$ is a prescribed involution.

Since $h$ is an involution, the pairs $\{x,x+h\}$ form a perfect
matching of the elements of $\Gamma$.
For an $h$-compatible partition $\mathcal P$, define
$Q_h(\mathcal P)$ to be the multigraph whose vertices are the
parts of $\mathcal P$ and whose edges are the pairs
$\{x,x+h\}$. Each pair contributes one edge joining the parts
containing its two elements. If every part has size $r$, then
$Q_h(\mathcal P)$ is loopless and $r$-regular. Parallel edges
are regarded as distinct edges.

For odd block sizes and a prescribed involution, we obtain the
following characterization.

\begin{theorem}\label{thm:partition}
Let $r\ge3$ be odd, let $\Gamma$ be a finite Abelian group, and
let $h\in I(\Gamma)$. There exists an $h$-compatible zero-sum
partition of $\Gamma$ of block size $r$ if and only if
\[
 r\mid|\Gamma|\qquad\text{and}\qquad |I(\Gamma)|\ge3.
\]
Moreover, when these conditions hold, the partition $\mathcal P$
can be chosen so that $Q_h(\mathcal P)$ is Hamiltonian and has
a proper edge coloring with $r$ colors.
\end{theorem}

The theorem applies to every involution $h$, including those
for which $\langle h\rangle$ is not a direct summand of $\Gamma$.
To prove Theorem~\ref{thm:partition}, we first construct, by
induction, a Hamiltonian multigraph with a bijective zero-sum
edge labeling by $\Gamma/\langle h\rangle$. We then assign the two elements
of each edge label to its endpoints and modify these assignments
so that the labels at every vertex sum to zero. The sets of
labels at the vertices form the required partition.

One application of our main theorem concerns complete mappings
of finite Abelian groups.
A \emph{complete mapping} of $\Gamma$ is a bijection
$\varphi\colon\Gamma\to\Gamma$ such that
$x\mapsto x+\varphi(x)$ is also bijective. A finite Abelian
group admits a complete mapping if and only if
$|I(\Gamma)|\ne1$~\cite{HallPaige1955}.
For results on complete mappings, see
\cite{Cichacz2025,Evans,Tan3,ref_MuyPok}.

A related problem concerns the cycle structure of orthomorphisms.
An \emph{orthomorphism} of $\Gamma$ is a bijection
$\theta\colon\Gamma\to\Gamma$ such that
$x\mapsto\theta(x)-x$ is also bijective. Every complete mapping
$\varphi$ gives an orthomorphism $\theta(x)=x+\varphi(x)$.
Friedlander, Gordon and Tannenbaum~\cite{Tan3} conjectured that,
for every finite Abelian group $\Gamma$ with $|I(\Gamma)|\ne1$
and every $k\ge2$ dividing $|\Gamma|-1$, there exists an
orthomorphism fixing $0$ and permuting the remaining elements
as a product of disjoint $k$-cycles.
M\"uyesser~\cite{Muyesser2026} proved this conjecture for all
sufficiently large groups.

For block size three, our partitions give complete mappings
directly. Every $h$-compatible zero-sum partition of $\Gamma$
into triples induces a complete mapping $\varphi$ whose cycles
are precisely the parts of the partition and which satisfies
$\varphi(x)\ne x+h$ for every $x\in\Gamma$.
We give the corresponding characterization in
Proposition~\ref{prop:triples}.

As another application of Theorem~\ref{thm:partition}, we
characterize the finite Abelian groups $\Gamma$ with
$r\mid|\Gamma|$ that admit an $r$-regular
$\Gamma$-distance-magic graph.
Such a graph exists if and only if $|I(\Gamma)|\ge3$.
Our construction makes each part of the partition a complete
graph and uses the pairs $\{x,x+h\}$ as a perfect matching
between the parts.
The resulting graph is Hamiltonian, $r$-regular, and admits
an $r$-partition.
For the same graph, $0$ and every involution of $\Gamma$
occur as magic constants. See Corollary~\ref{cor:magic}.

In Section~\ref{sec:edge}, we construct the zero-sum edge
labelings used in the proof of Theorem~\ref{thm:partition}.
We prove the theorem and discuss complete mappings for block
size three in Section~\ref{sec:partition}.
Section~\ref{sec:magic} contains the application to group
distance magic labeling and two examples.

\Needspace{12\baselineskip}
\section{Zero-sum edge labelings}\label{sec:edge}

We first construct a Hamiltonian multigraph whose edge labels sum
to zero at every vertex. This is the auxiliary construction used in the proof of
Theorem~\ref{thm:partition}.

For a loopless multigraph $Q$ and an Abelian group $A$, a
\emph{bijective zero-sum edge labeling} is a bijection
$f\colon E(Q)\to A$ such that
\begin{equation}\label{eq:edgezero}
 \sum_{e\ni v}f(e)=0
\end{equation}
for all $v\in V(Q)$.

\begin{proposition}\label{prop:edge}
Let $r\ge3$ be odd, and let $A$ be a finite Abelian group with
$r\mid |A|$. There exists a loopless $r$-regular multigraph $Q$
with a bijective zero-sum edge labeling by $A$ and an edge
decomposition
\begin{equation}\label{eq:decomposition}
 E(Q)=E(C)\,\dot\cup\,M_1\,\dot\cup\cdots\dot\cup\,M_{r-2},
\end{equation}
where $C$ is a Hamiltonian cycle and each $M_j$ is a perfect matching.
\end{proposition}

\begin{proof}
We prove the proposition by induction on $|A|$, with $r$ fixed.

Suppose first that $|A|=r$. Let $Q$ consist of two vertices joined
by $r$ parallel edges, and label these edges bijectively by the
elements of $A$. Both vertices have degree $r$, and the sum of
the labels incident with either vertex is
\[
 \sum_{a\in A}a=0.
\]
Indeed, since $|A|$ is odd, the nonzero elements of $A$ can be
partitioned into pairs $\{a,-a\}$, each of sum zero.
Any two edges form a Hamiltonian cycle $C$ of length two.
Each of the remaining $r-2$ edges is a perfect matching.
Thus~\eqref{eq:decomposition} holds.

Now suppose that $|A|>r$ and that the assertion holds for all
smaller groups whose orders are divisible by $r$.
Choose a prime $p$ dividing $|A|/r$, and let
$K=\langle k\rangle$ be a subgroup of $A$ of order $p$.
Then $|A/K|<|A|$ and $r$ divides $|A/K|$.

By the induction hypothesis, there exists a loopless $r$-regular
multigraph $H$ with a bijective zero-sum edge labeling
\[
 \bar f\colon E(H)\longrightarrow A/K
\]
and an edge decomposition
\[
 E(H)=E(C)\,\dot\cup\,M_1\,\dot\cup\cdots\dot\cup\,M_{r-2},
\]
where
\[
 C=v_0v_1\cdots v_{n-1}v_0
\]
is a Hamiltonian cycle and each $M_j$ is a perfect matching.
Since $H$ is $r$-regular and $r$ is odd, $n$ is even.

Let $e_i$ be the edge of $C$ joining $v_i$ to $v_{i+1}$,
and let $m_{i,j}$ be the edge of $M_j$ incident with $v_i$.
The subscripts of vertices and cycle edges are taken modulo $n$.
If $n=2$, then $e_0$ and $e_1$ are distinct parallel edges.

Put $\bar a_i=\bar f(e_i)$ and choose $a_i\in A$ such that
$a_i+K=\bar a_i$. Since the labels of the edges incident with
each vertex of $H$ sum to zero, we have
\begin{equation}\label{eq:quotientweights}
 \bar a_{i-1}+\bar a_i+
 \sum_{j=1}^{r-2}\bar f(m_{i,j})=0
\end{equation}
for all $i\in\mathbb Z_n$.
We consider two cases.

\medskip
\noindent\textbf{Case 1. $p$ is odd.}

Let $V(Q)=\{(i,t):i\in\mathbb Z_n,\ t\in\mathbb Z_p\}$.
First, construct a Hamiltonian cycle $C'$ by placing these
vertices in the cyclic order
\[
 (0,0),(1,0),\ldots,(n-1,0),
 (0,1),(1,1),\ldots,(n-1,1),
 \ldots,
 (0,p-1),\ldots,(n-1,p-1).
\]
Label the edge leaving $(i,t)$ in this order by $a_i+tk$.
For each fixed $i$, these labels use every element of $a_i+K$
exactly once.

We next construct the perfect matchings
$M'_2,\ldots,M'_{r-2}$.
For each $j\in\{2,\ldots,r-2\}$ and each edge
$e=v_iv_\ell\in M_j$, choose a representative $b_e\in A$ of
the coset $\bar f(e)$. Add the edges
\[
 (i,t)(\ell,t)
\]
for all $t\in\mathbb Z_p$, and label them by
\[
 b_e+(-1)^jtk,
\]
respectively. For each fixed $j$, these edges form a perfect
matching $M'_j$. Moreover, the edges corresponding to $e$ use
every element of $b_e+K$ exactly once.

Since $r$ is odd, the terms involving $tk$ cancel in pairs.
Thus the sum of the matching-edge labels at $(i,t)$ is
\[
 \sum_{j=2}^{r-2}
 \bigl(b_{m_{i,j}}+(-1)^jtk\bigr)
 =
 \sum_{j=2}^{r-2}b_{m_{i,j}}.
\]
Put
\[
 \sigma_i=\sum_{j=2}^{r-2}b_{m_{i,j}}.
\]
This sum does not depend on $t$. If $r=3$, no matchings are added
in this step, and we put $\sigma_i=0$.

Every vertex now has degree $r-1$. We add a perfect matching
$M'_1$ and choose its labels so that the sum at every vertex is zero.
We first determine the label needed at each vertex.
For $1\le i<n$, the two cycle edges incident with $(i,t)$ have
labels
\[
 a_{i-1}+tk
 \quad\text{and}\quad
 a_i+tk.
\]
Hence the sum of the labels already incident with $(i,t)$ is
\[
 a_{i-1}+a_i+\sigma_i+2tk.
\]
At $(0,t)$, the incoming cycle edge comes from $(n-1,t-1)$,
so the corresponding sum is
\[
 a_{n-1}+a_0+\sigma_0+(2t-1)k.
\]
Here and below, the second coordinates are taken modulo $p$.
Thus the label required on the last edge at $(i,t)$ is
\[
 \rho_i(t)=
 \begin{cases}
 -a_{i-1}-a_i-\sigma_i-2tk,
       &1\le i<n,\\[2pt]
 -a_{n-1}-a_0-\sigma_0-(2t-1)k,
       &i=0.
 \end{cases}
\]

We now construct $M'_1$.
By~\eqref{eq:quotientweights}, for every $i$ and $t$,
\[
 \rho_i(t)+K
 =
 -\bar a_{i-1}-\bar a_i
 -\sum_{j=2}^{r-2}\bar f(m_{i,j})
 =
 \bar f(m_{i,1}).
\]
Thus every $\rho_i(t)$ belongs to the coset $\bar f(m_{i,1})$.
Moreover, if
$\rho_i(t)=\rho_i(t')$, then
\[
 2(t-t')k=0.
\]
Since $k$ has odd prime order $p$, this implies $t=t'$ in
$\mathbb Z_p$. Thus, for each fixed $i$, the values
$\rho_i(t)$ run through the coset $\bar f(m_{i,1})$ exactly once.

Consider an edge $e=v_iv_\ell\in M_1$.
Since $m_{i,1}=m_{\ell,1}=e$, both lists
\[
 \{\rho_i(t):t\in\mathbb Z_p\}
 \quad\text{and}\quad
 \{\rho_\ell(s):s\in\mathbb Z_p\}
\]
contain every element of $\bar f(e)$ exactly once.
For each $z\in\bar f(e)$, let $t$ and $s$ be the unique indices
such that
\[
 \rho_i(t)=z=\rho_\ell(s).
\]
Add an edge joining $(i,t)$ to $(\ell,s)$ and label it by $z$.
The sum of the incident edge labels is now zero at both endpoints.

Doing this for every edge of $M_1$ gives a perfect matching
$M'_1$. Each vertex $(i,t)$ belongs to exactly one of the new
edges. Moreover, these edges are not loops, since $H$ is loopless
and hence $i\ne\ell$.

The labeling is bijective. Indeed, the edges of $C'$ use all elements
of the cosets assigned to the edges of $C$. Similarly, for each $j$,
the edges of $M'_j$ use all elements of the cosets assigned to the
edges of $M_j$. Each element of each coset is used exactly once.
Since $\bar f$ is bijective, these cosets partition $A$.

The multigraph $Q$ is loopless and $r$-regular, the labels of the
edges incident with every vertex sum to zero, and
\[
 E(Q)=E(C')\,\dot\cup\,M'_1\,\dot\cup\cdots\dot\cup\,M'_{r-2}.
\]
This completes Case~1.

\medskip
\noindent\textbf{Case 2. $p=2$.}

In this case $K=\{0,k\}$ and $2k=0$.
For each $1\le j\le r-2$ and each matching edge
$e=v_iv_\ell\in M_j$, assign the two elements of the coset
$\bar f(e)$ to its endpoints, writing them as
$b_{i,j}$ and $b_{\ell,j}$, where $b_{\ell,j}=b_{i,j}+k$.
Thus every element of every coset assigned to a matching edge
occurs exactly once among the $b_{i,j}$.

For each $i$, equation~\eqref{eq:quotientweights} implies that
\[
 a_{i-1}+a_i+\sum_{j=1}^{r-2}b_{i,j}\in K.
\]
Hence there is a unique $c_i\in\mathbb Z_2$ such that
\begin{equation}\label{eq:binarydefect}
 a_{i-1}+a_i+\sum_{j=1}^{r-2}b_{i,j}+c_i k=0.
\end{equation}

Replace each vertex $v_i$ by two vertices $(i,0)$ and $(i,1)$.
For each $i\in\mathbb Z_n$ and $t\in\mathbb Z_2$, add an edge
\begin{equation}\label{eq:binarycycle}
 (i,t)(i+1,t+c_i)
\end{equation}
with the label $a_i+(t+c_i)k$.
For each $i$ and $j\in\{1,\ldots,r-2\}$, also add an edge
$f_{i,j}$ joining $(i,0)$ to $(i,1)$, and label it by $b_{i,j}$.
The edges $f_{i,1},\ldots,f_{i,r-2}$ are distinct parallel edges.

Each vertex $(i,t)$ is incident with two edges of the form
\eqref{eq:binarycycle} and the $r-2$ edges $f_{i,j}$.
Thus the resulting multigraph $Q$ is $r$-regular.
It is loopless because consecutive vertices of $C$ are distinct,
and $(i,0)\ne(i,1)$.

The two edges corresponding to $e_i$ receive the two distinct
labels in $a_i+K$. The edges $f_{i,j}$ receive all elements
of the cosets assigned to the old matching edges, each once.
Since the old edge labels are distinct cosets, the new labeling
is bijective.

The incoming edge of the form~\eqref{eq:binarycycle} at $(i,t)$
has label $a_{i-1}+tk$, and the outgoing edge has label
$a_i+(t+c_i)k$. Therefore, the incident label sum at $(i,t)$ is
\[
 \begin{aligned}
 &(a_{i-1}+tk)+(a_i+(t+c_i)k)
       +\sum_{j=1}^{r-2}b_{i,j}\\
 &\qquad =
 a_{i-1}+a_i+\sum_{j=1}^{r-2}b_{i,j}+c_i k
 =0,
 \end{aligned}
\]
where we used $2k=0$ and~\eqref{eq:binarydefect}.

It remains to find the Hamiltonian cycle and the perfect matchings.
Let $D$ be the spanning subgraph consisting of the edges
in~\eqref{eq:binarycycle}, and put
\[
 F_j=\{f_{i,j}:i\in\mathbb Z_n\}
\]
for each $1\le j\le r-2$.
Every vertex has degree two in $D$, and each $F_j$ is a perfect
matching.

Following the edges of $D$ once around the old cycle changes
the second coordinate by
\[
 c_0+c_1+\cdots+c_{n-1}\in\mathbb Z_2.
\]
If this sum is $1$, then after $n$ steps we reach the other copy
of the starting vertex, and after another $n$ steps we return
to the starting vertex. Thus $D$ is a Hamiltonian cycle, and
\[
 E(Q)=E(D)\,\dot\cup\,F_1\,\dot\cup\cdots\dot\cup\,F_{r-2}
\]
is the required decomposition.

Suppose instead that the sum is $0$. Then $D$ consists of two
vertex-disjoint cycles of length $n$, each containing exactly
one of $(i,0)$ and $(i,1)$ for every $i$.

Choose an edge $e_i$ of the old cycle $C$ and remove from $D$
the two edges corresponding to $e_i$. This turns the two cycles
into two vertex-disjoint paths. Each path has one endpoint
whose first coordinate is $i$ and one whose first coordinate
is $i+1$.
The edge $f_{i,1}$ joins the two endpoints with first coordinate
$i$, and $f_{i+1,1}$ joins the two endpoints with first coordinate
$i+1$. Adding these two edges to the paths gives a Hamiltonian
cycle $C'$.

The spanning multigraph with edge set $E(D)\cup F_1$ is
$3$-regular, whereas $C'$ has degree two at every vertex.
Consequently,
\[
 M'_1=\bigl(E(D)\cup F_1\bigr)\setminus E(C')
\]
is a perfect matching. Hence
\[
 E(Q)=E(C')\,\dot\cup\,M'_1\,\dot\cup\,
 F_2\,\dot\cup\cdots\dot\cup\,F_{r-2}
\]
is the required decomposition.
This completes the induction.

\end{proof}

\section{Compatible zero-sum partitions}\label{sec:partition}

We now use Proposition~\ref{prop:edge} to prove our main theorem.
The edge labels are elements of $\Gamma/\langle h\rangle$,
and the elements of $\Gamma$ are assigned to the edge ends.

\begin{proof}[Proof of Theorem~\ref{thm:partition}]
Suppose first that such a partition $\mathcal P$ exists.
Since every part has $r$ elements, $r$ divides $|\Gamma|$.
Moreover, the parts are disjoint, cover $\Gamma$, and have sum zero.
Therefore
\[
 \sum_{x\in\Gamma}x
 =\sum_{T\in\mathcal P}\sum_{x\in T}x=0.
\]
By~Lemma~\ref{CN}, $\Gamma$ cannot have exactly one involution.
Since $h\in I(\Gamma)$, the group has at least one involution,
and hence $|I(\Gamma)|\ge3$.

Conversely, suppose that $r$ divides $|\Gamma|$ and $|I(\Gamma)|\ge3$.
The subgroup $I(\Gamma)\cup\{0\}$ has order at least four, so
$4$ divides $|\Gamma|$. Since $r$ is odd, it follows that 
$4r$ divides $|\Gamma|$. Also, by Lemma~\ref{CN}, $ \sum_{x\in\Gamma}x=0.$ Put
\[
 H=\langle h\rangle=\{0,h\},\qquad A=\Gamma/H,
\]
and let $\pi\colon\Gamma\to A$ be the quotient map.
Then $|A|=|\Gamma|/2$, so $r$ divides $|A|$ and $|A|\ge2r$.
By Proposition~\ref{prop:edge}, there exists a loopless
$r$-regular multigraph $Q$ with a bijective zero-sum edge labeling
\[
 f\colon E(Q)\longrightarrow A
\]
and an edge decomposition
\[
 E(Q)=E(C)\,\dot\cup\,M_1\,\dot\cup\cdots\dot\cup\,M_{r-2},
\]
where $C$ is a Hamiltonian cycle and each $M_j$ is a perfect matching.

We now label the two ends of every edge of $Q$.
For an edge $e=uv$, the label $f(e)$ is a coset of $H$ with two
elements, say $x$ and $x+h$. Assign one element to the end at $u$
and the other to the end at $v$. Denote the label at the $v$-end
of $e$ by $\lambda(v,e)$. Thus
\[
 \{\lambda(u,e),\lambda(v,e)\}=\{x,x+h\},
 \qquad
 \pi(\lambda(u,e))=\pi(\lambda(v,e))=f(e).
\]
Since $f$ is bijective, every element of $\Gamma$ is used exactly
once. In particular, the $r$ edge ends at each vertex have distinct
labels.

For $v\in V(Q)$, let
\[
 s(v)=\sum_{e\ni v}\lambda(v,e).
\]
Since
\[
 \pi(s(v))
 =\sum_{e\ni v}\pi(\lambda(v,e))
 =\sum_{e\ni v}f(e)=0,
\]
we have $s(v)\in\ker\pi=H=\{0,h\}$ for every $v$.
Also, since every element of $\Gamma$ occurs at exactly one edge end,
\begin{equation}\label{eq:defecttotal}
 \sum_{v\in V(Q)}s(v)=\sum_{x\in\Gamma}x=0.
\end{equation}

We now modify the assignments so that $s(v)=0$ for every vertex.
Consider an edge $e=uv$ whose end labels are $x$ and $x+h$.
Interchanging these labels adds $h$ to the sum at $u$ and also
to the sum at $v$, since $-h=h$.
Thus the interchange changes $0$ to $h$ and $h$ to $0$ at both
endpoints. Every group element is still used once, so
\eqref{eq:defecttotal} still holds.

Delete one edge of $C$ to obtain a Hamiltonian path
$v_1v_2\cdots v_q$. For $i=1,\ldots,q-1$, proceed as follows.
If $s(v_i)=0$, make no change. If $s(v_i)=h$, interchange the
end labels of the edge $v_iv_{i+1}$.
After this step, $s(v_i)=0$. The interchange affects only $v_i$
and $v_{i+1}$, so the sums at $v_1,\ldots,v_{i-1}$ remain zero.
Consequently, after the last step,
\[
 s(v_1)=\cdots=s(v_{q-1})=0.
\]
Equation~\eqref{eq:defecttotal} then gives $s(v_q)=0$ as well.

After these interchanges, put
\[
 T_v=\{\lambda(v,e):e\ni v\}
\]
for all $v\in V(Q)$.
The graph $Q$ is $r$-regular, and all assigned elements are distinct,
so $|T_v|=r$. Moreover,
\[
 \sum_{x\in T_v}x=s(v)=0.
\]
Since every element of $\Gamma$ is assigned to exactly one edge end,
the sets $T_v$ are pairwise disjoint and cover $\Gamma$.
Therefore
\[
 \mathcal P=\{T_v:v\in V(Q)\}
\]
is a zero-sum partition of block size $r$.

To verify $h$-compatibility, let $x\in T_v$, and let $e$ be the
edge whose $v$-end is assigned $x$. The other end of $e$ is
assigned $x+h$. Since $Q$ is loopless, that end belongs to a
vertex $u\ne v$. Hence $x+h\in T_u$ and $x+h\notin T_v$.
It follows that $T_v\cap(T_v+h)=\varnothing$ for every $v$.

Finally, the map $v\mapsto T_v$ is an isomorphism from $Q$ to
$Q_h(\mathcal P)$. Indeed, an edge $e=uv$ corresponds to the pair
$\{\lambda(u,e),\lambda(v,e)\}=\{x,x+h\}$ joining $T_u$ to $T_v$.
Conversely, every such pair comes from exactly one edge of $Q$.
Thus $Q_h(\mathcal P)$ is Hamiltonian and hence connected.

Since $Q$ is $r$-regular and $r$ is odd, $q$ is even.
Color the edges of $C$ alternately with colors $1$ and $2$,
and give every edge of $M_j$ color $j+2$ for
$1\le j\le r-2$. At each vertex, the two cycle edges have different
colors, and there is exactly one incident edge from each $M_j$.
This gives a proper edge coloring of $Q$, and hence of
$Q_h(\mathcal P)$, with $r$ colors.
\end{proof}

For block size three, compatible zero-sum partitions can be
described by complete mappings with an additional condition.
Here $\varphi^2=\varphi\circ\varphi$.

\begin{proposition}\label{prop:triples}
Let $\Gamma$ be a finite Abelian group and let $h\in I(\Gamma)$.
There exists an $h$-compatible zero-sum partition of $\Gamma$
into triples if and only if there is a bijection
$\varphi\colon\Gamma\to\Gamma$ such that, for every $x\in\Gamma$,
\begin{equation}\label{eq:triple-conditions}
 \varphi(x)\ne x,\qquad
 x+\varphi(x)+\varphi^2(x)=0,\qquad
 \varphi(x)\ne x+h.
\end{equation}
Every such $\varphi$ is a complete mapping, and its cycles have
length three and form an $h$-compatible zero-sum partition.
\end{proposition}

\begin{proof}
Let $\mathcal P$ be an $h$-compatible zero-sum partition into
triples. Cyclically permute the elements of each part to obtain
a bijection $\varphi$. Its cycles have length three, so it has
no fixed points. The zero-sum condition gives
$x+\varphi(x)+\varphi^2(x)=0$, and $h$-compatibility gives
$\varphi(x)\ne x+h$.

Conversely, suppose that a bijection $\varphi$ satisfies all
three conditions in~\eqref{eq:triple-conditions}.
Applying the middle identity in~\eqref{eq:triple-conditions}
at $x$ and at $\varphi(x)$ gives
\[
 x+\varphi(x)+\varphi^2(x)=0
 \quad\text{and}\quad
 \varphi(x)+\varphi^2(x)+\varphi^3(x)=0.
\]
Subtracting yields $\varphi^3(x)=x$.
Since $\varphi$ has no fixed points, every cycle has length
three, and its elements sum to zero.
If $\varphi^2(x)=x+h$, then
$\varphi(x+h)=\varphi^3(x)=x=(x+h)+h$, contrary to the last
condition in~\eqref{eq:triple-conditions} at $x+h$.
Thus neither $\varphi(x)$ nor $\varphi^2(x)$ equals $x+h$.
The cycles therefore form an $h$-compatible zero-sum partition.

Finally, every bijection $\varphi$ satisfying
\eqref{eq:triple-conditions} satisfies
\[
 x+\varphi(x)=-\varphi^2(x).
\]
Since the right-hand side is a bijection of $\Gamma$,
$\varphi$ is a complete mapping.
\end{proof}

\section{Application to group distance magic labeling}\label{sec:magic}

We apply Theorem~\ref{thm:partition} to construct Hamiltonian
regular $\Gamma$-distance-magic graphs.
For a graph $G$ and a vertex $x\in V(G)$, let $N_G(x)$ denote
the set of neighbors of $x$.

Froncek~\cite{Fro} introduced the notion of a group distance magic
labeling. Let $(\Gamma,+)$ be an Abelian group of order $|V(G)|$.
A \emph{$\Gamma$-distance-magic labeling} of $G$ is a bijection
$\ell\colon V(G)\to\Gamma$ for which there exists $\mu\in\Gamma$
such that
\[
 \sum_{y\in N_G(x)}\ell(y)=\mu
\]
for all $x\in V(G)$.
The element $\mu$ is called the \emph{magic constant} of the labeling.
If $G$ admits such a labeling, then $G$ is called
\emph{$\Gamma$-distance-magic}.

Cichacz and Froncek~\cite{CicFro} proved the following necessary
condition for regular graphs of odd degree.

\begin{theorem}[\cite{CicFro}]\label{gr:odd}
Let $G$ be an $r$-regular graph of order $n$, where $r$ is odd.
If $\Gamma$ is an Abelian group of order $n$ with
$|I(\Gamma)|=1$, then $G$ is not $\Gamma$-distance-magic.
\end{theorem}

Since $r$ is odd, every $r$-regular graph has even order. Hence,
if $G$ is $\Gamma$-distance-magic, then $\Gamma$ has an involution.
By Theorem~\ref{gr:odd}, it cannot have exactly one involution.
Thus $|I(\Gamma)|\ge3$.

For a positive integer $t$, let $tG$ denote the disjoint union of
$t$ copies of $G$. If $r$ divides $|\Gamma|$ and $|I(\Gamma)|\ge3$,
Cichacz~\cite[Theorem~15]{Cichacz2018} showed that
$\frac{|\Gamma|}{2r}K_{r,r}$ is $\Gamma$-distance-magic.
This gives a disconnected $r$-regular $\Gamma$-distance-magic graph.
In the cubic case, Cichacz and Miklavi\v{c}~\cite{CichaczMiklavic2026}
obtained nonexistence results for generalized Petersen graphs and
gave constructions for disjoint unions of group distance magic
graphs. Odd regular $(\mathbb{Z}_2)^m$-distance magic graphs were recently studied
in~\cite{CicGroZua}.

We use the following definition of Cichacz and
Nikodem~\cite{CicNik}. Let $p\ge2$. A partition
$V(G)=V_1\,\dot\cup\cdots\dot\cup\,V_p$ into nonempty sets is
called a \emph{$p$-partition} of $G$ if
\[
 |N_G(x)\cap V_1|=\cdots=|N_G(x)\cap V_p|
\]
for all $x\in V(G)$.
If $G$ is $r$-regular and has a $p$-partition, then $p$ divides $r$,
and every vertex has exactly $r/p$ neighbors in each part.

Let $\mathcal P$ be an $h$-compatible partition of $\Gamma$ into
$r$-element sets. Define $G(\mathcal P,h)$ on vertex set $\Gamma$
by making each part a complete graph and adding the edges
$\{x,x+h\}$. These additional edges form a perfect matching.
Since $x$ and $x+h$ belong to different parts, the graph is simple
and $r$-regular. The complete graphs form a $K_r$-factor, that is,
a spanning union of vertex-disjoint copies of $K_r$.
Contracting each copy of $K_r$ to one vertex gives $Q_h(\mathcal P)$.

The zero-sum condition gives a group distance magic labeling,
while the Hamiltonian cycle and the edge coloring of
$Q_h(\mathcal P)$ give additional properties of $G(\mathcal P,h)$.
Together with Theorem~\ref{gr:odd}, this yields the following
characterization.

\begin{corollary}\label{cor:magic}
Let $r\ge3$ be odd, and let $\Gamma$ be a finite Abelian group
with $r\mid|\Gamma|$.
There exists an $r$-regular $\Gamma$-distance-magic graph if
and only if $|I(\Gamma)|\ge3$.
Moreover, when these conditions hold, there is a single
Hamiltonian $r$-regular graph with an $r$-partition that admits
a $\Gamma$-distance-magic labeling with magic constant $c$
for every $c\in\{0\}\cup I(\Gamma)$.
\end{corollary}

\begin{proof}
Suppose that an $r$-regular $\Gamma$-distance-magic graph exists.
Since $r$ is odd, its order $|\Gamma|$ is even, so $\Gamma$
has an involution. Theorem~\ref{gr:odd} excludes exactly one
involution, and hence $|I(\Gamma)|\ge3$.

Conversely, suppose that $|I(\Gamma)|\ge3$ and fix
$h\in I(\Gamma)$. By Theorem~\ref{thm:partition}, there
exists an $h$-compatible zero-sum partition $\mathcal P$ of block
size $r$ such that $Q_h(\mathcal P)$ is Hamiltonian and has a proper
edge coloring with $r$ colors.
Let $G=G(\mathcal P,h)$. For $x\in T\in\mathcal P$,
\[
 N_G(x)=(T\setminus\{x\})\cup\{x+h\},\qquad
 \sum_{y\in N_G(x)}y=-x+(x+h)=h.
\]
Thus the identity labeling is $\Gamma$-distance-magic with magic
constant $h$.

We next show that $G$ is Hamiltonian. Take a Hamiltonian cycle
of $Q_h(\mathcal P)$. At each part $T\in\mathcal P$, the two edges
of this cycle correspond to matching edges with distinct endpoints
in $T$. Since $T$ induces a complete graph, there is a path through
all vertices of $T$ joining these endpoints. These paths, together
with the matching edges corresponding to the cycle, form a
Hamiltonian cycle of $G$.

For each $c\in\{0\}\cup I(\Gamma)$, define
$\ell_c(x)=x+h+c$ on the same graph $G$. Translation is a bijection,
and, since $r$ is odd and $2h=2c=0$, we have
\[
 \sum_{y\in N_G(x)}\ell_c(y)=h+r(h+c)=c.
\]
It remains to find an $r$-partition of $G$.
Take a proper edge coloring of $Q_h(\mathcal P)$ with $r$ colors.
Give both vertices $x$ and $x+h$ the color of the edge corresponding
to $\{x,x+h\}$. Since $Q_h(\mathcal P)$ is $r$-regular, every part
$T$ contains one vertex of each color. Thus each vertex of $G$ has
one neighbor of every other color in its part and one neighbor of
its own color along the matching.
The color classes $V_1,\ldots,V_r$ satisfy
\[
 |N_{G(\mathcal P,h)}(x)\cap V_i|=1
\]
for all $x\in\Gamma$ and all $1\le i\le r$.
This is an $r$-partition of $G$.
\end{proof}

The $r$-partition in Corollary~\ref{cor:magic} also allows us to
apply the following result to disjoint copies of the constructed
graph.

\begin{lemma}[\cite{CichaczMiklavic2026}]\label{lemgl}
Let $A$ and $B$ be finite Abelian groups, and put $t=|B|$.
Suppose that $G$ is an $r$-regular $A$-distance-magic graph having
a $p$-partition, where $p\ge2$. If $p$ is even or $|I(B)|\ne1$,
then $tG$ is $(A\times B)$-distance-magic.
\end{lemma}

In particular, if $B$ is a finite Abelian group with
$|I(B)|\ne1$, then $|B|G$ is $(\Gamma\times B)$-distance-magic
for the graph $G$ in Corollary~\ref{cor:magic}.

Our construction requires $r\mid|\Gamma|$ because it uses a
$K_r$-factor. This divisibility is not necessary for every
$r$-regular $\Gamma$-distance-magic graph. For example,
$K_{3,3,3,3}$ is $9$-regular of order $12$ and admits a
$\mathbb Z_2^2\oplus\mathbb Z_3$-distance-magic labeling
\cite[Section~4]{Cichacz2018}.

We give two examples of the partitions in Theorem~\ref{thm:partition}
and their corresponding group distance magic graphs.
In the first, $\langle h\rangle$ is not a direct summand of $\Gamma$.

\begin{example}\label{ex:cube}
Let $r=3$, $\Gamma=\mathbb Z_2\oplus\mathbb Z_{12}$, and $h=(0,6)$.
The following triples form the $h$-compatible zero-sum partition
$\mathcal P=\{T_0,\ldots,T_7\}$.
\begin{align*}
 T_0&=\{(0,0),(0,8),(0,4)\},&
 T_1&=\{(0,6),(1,11),(1,7)\},\\
 T_2&=\{(0,2),(0,3),(0,7)\},&
 T_3&=\{(1,5),(0,9),(1,10)\},\\
 T_4&=\{(0,10),(1,0),(1,2)\},&
 T_5&=\{(1,1),(1,6),(0,5)\},\\
 T_6&=\{(0,1),(1,8),(1,3)\},&
 T_7&=\{(1,4),(0,11),(1,9)\}.
\end{align*}
In $Q_h(\mathcal P)$, the vertices $T_i$ and $T_j$ are adjacent
exactly when the three-digit binary expansions of $i$ and $j$ differ
in one position. Thus the quotient is the cube graph. Here
$h=2(0,3)\in I(\Gamma)$. The graph $G(\mathcal P,h)$ is a
Hamiltonian cubic $\Gamma$-distance-magic graph on $24$ vertices
with magic constant $(0,6)$.
The graph and its quotient are shown in
Figure~\ref{fig:cube-example}.
\end{example}

In the second example, $r=5$ and $|\Gamma|=20$, the smallest
possible group order in Theorem~\ref{thm:partition} for $r=5$.

\begin{example}\label{ex:five}
Write $\mathbb Z_2^2=\{0,a,b,c\}$, where $a+b=c$, and take
$\Gamma=\mathbb Z_2^2\oplus\mathbb Z_5$ and $h=(a,0)$.
An $h$-compatible zero-sum partition
$\mathcal P=\{T_0,T_a,T_b,T_c\}$ is given by
\begin{align*}
 T_0&=\{(0,0),(0,1),(0,2),(0,3),(0,4)\},\\
 T_a&=\{(a,0),(b,1),(c,2),(c,3),(c,4)\},\\
 T_b&=\{(b,0),(c,1),(a,2),(a,3),(a,4)\},\\
 T_c&=\{(c,0),(a,1),(b,2),(b,3),(b,4)\}.
\end{align*}
The quotient has the Hamiltonian cycle $T_0T_aT_bT_cT_0$.
The pairs $T_0,T_b$ and $T_a,T_c$ are each joined by three edges,
and every other pair is joined by one edge.
The multigraph $Q_h(\mathcal P)$ has four vertices and degree five,
so it has parallel edges. The graph $G(\mathcal P,h)$ is simple,
Hamiltonian and $5$-regular. The identity labeling is
$\Gamma$-distance-magic with magic constant $(a,0)$.
The graph and its quotient are shown in
Figure~\ref{fig:five-example}.
\end{example}

\begin{figure}[htbp]
\centering
\begin{minipage}{.80\linewidth}
\centering
\begin{minipage}[c]{.70\linewidth}
\centering
\begin{tikzpicture}[
  x=.56cm,y=.56cm,
  vertex/.style={circle,draw=black!70,fill=white,
    inner sep=.5pt,minimum size=6.5mm,font=\fontsize{7}{8}\selectfont},
  clique/.style={draw=black!60,line width=.55pt},
  matching/.style={draw=blue!70!black,line width=.9pt}
]
\coordinate (p00) at (-4.4722,5.5500);
\coordinate (p01) at (-5.5500,4.4722);
\coordinate (p02) at (-4.0777,4.0777);
\coordinate (p10) at (4.4722,5.5500);
\coordinate (p11) at (5.5500,4.4722);
\coordinate (p12) at (4.0777,4.0777);
\coordinate (p20) at (-5.5500,-4.4722);
\coordinate (p21) at (-4.4722,-5.5500);
\coordinate (p22) at (-4.0777,-4.0777);
\coordinate (p30) at (5.5500,-4.4722);
\coordinate (p31) at (4.4722,-5.5500);
\coordinate (p32) at (4.0777,-4.0777);
\coordinate (p40) at (-2.5223,2.5223);
\coordinate (p41) at (-1.0500,2.1278);
\coordinate (p42) at (-2.1278,1.0500);
\coordinate (p50) at (2.5223,2.5223);
\coordinate (p51) at (1.0500,2.1278);
\coordinate (p52) at (2.1278,1.0500);
\coordinate (p60) at (-2.5223,-2.5223);
\coordinate (p61) at (-2.1278,-1.0500);
\coordinate (p62) at (-1.0500,-2.1278);
\coordinate (p70) at (2.5223,-2.5223);
\coordinate (p71) at (2.1278,-1.0500);
\coordinate (p72) at (1.0500,-2.1278);
\foreach \t in {0,...,7}{
  \draw[clique] (p\t0) -- (p\t1) -- (p\t2) -- cycle;
}
\foreach \u/\v in {00/10,60/22,20/01,21/31,02/40,52/71,
                    41/51,50/12,42/61,62/72,70/32,30/11}{
  \draw[matching] (p\u) -- (p\v);
}
\foreach \id/\lab in {
  00/{(0,0)},01/{(0,8)},02/{(0,4)},
  10/{(0,6)},11/{(1,11)},12/{(1,7)},
  20/{(0,2)},21/{(0,3)},22/{(0,7)},
  30/{(1,5)},31/{(0,9)},32/{(1,10)},
  40/{(0,10)},41/{(1,0)},42/{(1,2)},
  50/{(1,1)},51/{(1,6)},52/{(0,5)},
  60/{(0,1)},61/{(1,8)},62/{(1,3)},
  70/{(1,4)},71/{(0,11)},72/{(1,9)}}{
  \node[vertex] at (p\id) {$\lab$};
}
\foreach \t/\x/\y in {
  0/-6.8/4.7,1/6.8/4.7,2/-6.8/-4.7,3/6.8/-4.7,
  4/-1.7/3.4,5/1.7/3.4,6/-1.7/-3.4,7/1.7/-3.4}{
  \node[font=\scriptsize] at (\x,\y) {$T_{\t}$};
}
\end{tikzpicture}
\par\smallskip
{\small (a) The graph $G(\mathcal P,h)$.}
\end{minipage}\hfill
\begin{minipage}[c]{.28\linewidth}
\centering
\begin{tikzpicture}[
  x=.27cm,y=.27cm,
  qvertex/.style={circle,draw=black!65,fill=white,
    minimum size=6mm,inner sep=1pt,font=\scriptsize},
  qedge/.style={draw=blue!70!black,line width=.7pt}
]
\node[qvertex] (t0) at (-4.7,4.7) {$T_0$};
\node[qvertex] (t1) at (4.7,4.7) {$T_1$};
\node[qvertex] (t2) at (-4.7,-4.7) {$T_2$};
\node[qvertex] (t3) at (4.7,-4.7) {$T_3$};
\node[qvertex] (t4) at (-1.9,1.9) {$T_4$};
\node[qvertex] (t5) at (1.9,1.9) {$T_5$};
\node[qvertex] (t6) at (-1.9,-1.9) {$T_6$};
\node[qvertex] (t7) at (1.9,-1.9) {$T_7$};
\foreach \u/\v in {0/1,0/2,0/4,1/3,1/5,2/3,2/6,
                    3/7,4/5,4/6,5/7,6/7}{
  \draw[qedge] (t\u) -- (t\v);
}
\end{tikzpicture}
\par\smallskip
{\small (b) The quotient\\$Q_h(\mathcal P)$.}
\end{minipage}
\caption{A cubic $\Gamma$-distance-magic graph $G(\mathcal P,h)$ on
$24$ vertices, where $\Gamma=\mathbb Z_2\oplus\mathbb Z_{12}$ and
$h=(0,6)$. The triangles $T_0,\ldots,T_7$ are the zero-sum blocks
of $\mathcal P$. Blue edges form the perfect matching
$\{\{x,x+h\}:x\in\Gamma\}$, and gray edges lie within the blocks.
The identity labeling has magic constant $h$.
The quotient $Q_h(\mathcal P)$ in panel~(b) is the cube graph.}
\label{fig:cube-example}
\end{minipage}
\end{figure}

\begin{figure}[htbp]
\centering
\begin{minipage}{.80\linewidth}
\centering
\begin{minipage}[c]{.70\linewidth}
\centering
\begin{tikzpicture}[
  x=.60cm,y=.60cm,
  vertex/.style={circle,draw=black!70,fill=white,
    inner sep=.5pt,minimum size=6.5mm,font=\fontsize{7}{8}\selectfont},
  clique/.style={draw=black!60,line width=.55pt},
  matching/.style={draw=blue!70!black,line width=.9pt}
]
\coordinate (p00) at (-4.55,4.45);
\coordinate (p01) at (-2.65,4.45);
\coordinate (p02) at (-4.90,2.70);
\coordinate (p03) at (-3.60,1.95);
\coordinate (p04) at (-2.30,2.70);
\coordinate (pa0) at (4.55,4.45);
\coordinate (pa1) at (2.65,4.45);
\coordinate (pa2) at (4.90,2.70);
\coordinate (pa3) at (3.60,1.95);
\coordinate (pa4) at (2.30,2.70);
\coordinate (pb0) at (-4.55,-4.45);
\coordinate (pb1) at (-2.65,-4.45);
\coordinate (pb2) at (-4.90,-2.70);
\coordinate (pb3) at (-3.60,-1.95);
\coordinate (pb4) at (-2.30,-2.70);
\coordinate (pc0) at (4.55,-4.45);
\coordinate (pc1) at (2.65,-4.45);
\coordinate (pc2) at (4.90,-2.70);
\coordinate (pc3) at (3.60,-1.95);
\coordinate (pc4) at (2.30,-2.70);
\foreach \t in {0,a,b,c}{
  \foreach \u/\v in {0/1,0/2,0/3,0/4,1/2,1/3,1/4,2/3,2/4,3/4}{
    \draw[clique] (p\t\u) -- (p\t\v);
  }
}
\draw[matching] (p00) .. controls (-4.55,6.35) and (4.55,6.35) .. (pa0);
\draw[matching] (p01) .. controls (-.50,5.60) and (.50,-5.60) .. (pc1);
\foreach \u in {2,3,4}{
  \draw[matching] (p0\u) -- (pb\u);
  \draw[matching] (pc\u) -- (pa\u);
}
\draw[matching] (pb0) .. controls (-4.55,-6.35) and (4.55,-6.35) .. (pc0);
\draw[matching] (pa1) .. controls (.50,5.60) and (-.50,-5.60) .. (pb1);
\foreach \id/\lab in {
  00/{(0,0)},01/{(0,1)},02/{(0,2)},03/{(0,3)},04/{(0,4)},
  a0/{(a,0)},a1/{(b,1)},a2/{(c,2)},a3/{(c,3)},a4/{(c,4)},
  b0/{(b,0)},b1/{(c,1)},b2/{(a,2)},b3/{(a,3)},b4/{(a,4)},
  c0/{(c,0)},c1/{(a,1)},c2/{(b,2)},c3/{(b,3)},c4/{(b,4)}}{
  \node[vertex] at (p\id) {$\lab$};
}
\node[font=\scriptsize] at (-5.6,3.5) {$T_0$};
\node[font=\scriptsize] at (5.6,3.5) {$T_a$};
\node[font=\scriptsize] at (-5.6,-3.5) {$T_b$};
\node[font=\scriptsize] at (5.6,-3.5) {$T_c$};
\end{tikzpicture}
\par\smallskip
{\small (a) The graph $G(\mathcal P,h)$.}
\end{minipage}\hfill
\begin{minipage}[c]{.28\linewidth}
\centering
\begin{tikzpicture}[
  x=.33cm,y=.33cm,
  qvertex/.style={circle,draw=black!65,fill=white,
    minimum size=6mm,inner sep=1pt,font=\scriptsize},
  qedge/.style={draw=blue!70!black,line width=.7pt}
]
\node[qvertex] (t0) at (-3.6,3.5) {$T_0$};
\node[qvertex] (ta) at (3.6,3.5) {$T_a$};
\node[qvertex] (tb) at (-3.6,-3.5) {$T_b$};
\node[qvertex] (tc) at (3.6,-3.5) {$T_c$};
\draw[qedge] (t0) -- (ta);
\draw[qedge] (t0) -- (tb);
\draw[qedge] (t0) -- (tc);
\draw[qedge] (ta) -- (tb);
\draw[qedge] (ta) -- (tc);
\draw[qedge] (tb) -- (tc);
\draw[qedge] (t0) to[bend left=18] (tb);
\draw[qedge] (t0) to[bend right=18] (tb);
\draw[qedge] (ta) to[bend left=18] (tc);
\draw[qedge] (ta) to[bend right=18] (tc);
\end{tikzpicture}
\par\smallskip
{\small (b) The quotient\\$Q_h(\mathcal P)$.}
\end{minipage}
\caption{A $5$-regular $\Gamma$-distance-magic graph $G(\mathcal P,h)$ on
$20$ vertices, where $\Gamma=\mathbb Z_2^2\oplus\mathbb Z_5$,
$\mathbb Z_2^2=\{0,a,b,c\}$ with $a+b=c$, and $h=(a,0)$.
Each block $T_0,T_a,T_b,T_c$ induces a copy of $K_5$.
Blue edges form the perfect matching $\{\{x,x+h\}:x\in\Gamma\}$,
and gray edges lie within the blocks.
The identity labeling has magic constant $h$.
In the quotient in panel~(b), the pairs $T_0,T_b$ and $T_a,T_c$
are each joined by three edges, and every other pair by one edge.}
\label{fig:five-example}
\end{minipage}
\end{figure}

\textbf{Generative AI statement.} During the preparation of this manuscript, the authors used ChatGPT (OpenAI) to assist with mathematical reasoning, including drafting some proofs, as well as language editing and the preparation of figures. All mathematical arguments were independently checked by the authors. The authors take full responsibility for the mathematical results, proofs, and all other content of the manuscript.
\FloatBarrier

\end{document}